\documentstyle[12pt]{article}
\newcommand{\sect}[1]{\section{#1}\setcounter{equation}{0}}

\font\mbn=msbm10 scaled \magstep1
\font\mbs=msbm7 scaled \magstep1
\font\mbss=msbm5 scaled \magstep1
\newfam\mbff
\textfont\mbff=\mbn
\scriptfont\mbff=\mbs
\scriptscriptfont\mbff=\mbss\def\mbf{\fam\mbff}
\def\Re{{\mbf R}}

\def\Co{{\mbf C}}

\newtheorem{Th}{Theorem}[section]
\newtheorem{Lm}[Th]{Lemma}
\newtheorem{C}[Th]{Corollary}

\newtheorem{R}[Th]{Remark}

\author{Alexander Brudnyi\thanks{Research supported in part by NSERC.
}\\
Department of Mathematics and Statistics\\
University of Calgary, Calgary\\
Canada}
\title{Hartogs Type Theorems on Coverings of Stein Manifolds}
\date{} 
\begin{document} 
\maketitle
\begin{abstract}
{We prove an analog of the classical Hartogs extension theorem for
certain (possibly unbounded) domains on coverings of Stein manifolds.
}
\end{abstract}

\sect{\hspace*{-1em}. Introduction.}
Let $D\subset\Co^{n}\ (n>1)$ be a bounded open set with a connected smooth 
boundary
$bD$. The classical Hartogs theorem states that any holomorphic function
in some neighbourhood of $bD$ can be extended to a holomorphic function on
a neighbourhood of the closure $\overline{D}$. The first rigorous proof of 
this theorem was given by 
Brown in 1936 
see, e.g., 
[F]. In [Bo] Bochner proved a similar extension result for functions defined 
on the $bD$ only. In modern language his result says that for a smooth 
function
defined on the $bD$ and satisfying the tangential Cauchy-Riemann equations
there is an extension to a holomorphic function in $D$ which is smooth
on $\overline{D}$. In fact, this statement follows from 
Bochner's proof (under some smoothness conditions). However at that time
there was not yet the notion of a $CR$-function. Over the years significant
contributions to the area of Hartogs theorem were made by many prominent
mathematicians, see the history and the references in the paper of
Harvey and Lawson [HL, Section 5]. A general Hartogs-Bochner type
theorem for bounded domains $D$ in Stein manifolds is proved by Harvey and
Lawson [HL, Theorem 5.1]. The proof relies heavily upon the fact that
for $n\geq 2$ any $\overline{\partial}$-equation with compact support on
a Stein manifold has a compactly supported solution.
In this paper we present a Hartogs type theorem for certain 
(possibly unbounded) domains on coverings of Stein manifolds which gives
an extension of the above cited result of [HL]. 
In order to formulate this theorem
we first introduce some notation and basic definitions.

Let $r: M'\to M$ be an unbranched covering of a Stein manifold $M$ of
complex dimension $n\geq 2$. Let
$D\subset M'$ be a domain with a connected piecewise smooth boundary $bD$
such that $r(D)\subset\subset M$. Assume that $M$ is equipped with a
Riemannian metric $g_{M}$. By $d$ we denote the path metric on $M'$ 
induced by the pullback of $g_{M}$. For a fixed $o\in D$ we set
$d_{o}(z):=d(o,z)$, $z\in M'$. Also, 
by $\overline{D}\subset M'$  we denote
the closure of $D$ and by ${\cal O}(D)$ the space of holomorphic functions
on $D$. Next, recall that a continuous
function $f$ on $bD$ is called $CR$ if for every smooth $(n,n-2)$-form
$\omega$ on $M'$ with a compact support one has
$$
\int_{bD}f\cdot\overline{\partial}\omega=0\ \! .
$$
If $f$ and $bD$ are smooth this is equivalent to $f$
being a solution of the tangential $CR$-equations: 
$\overline{\partial}_{b}f=0$ (see, e.g., [KR]).

Suppose that $f\in C(bD)$ is a $CR$-function satisfying for
some positive numbers $c_{1},c_{2}, \delta$ the following conditions
\begin{itemize}
\item[(1)]
$$
|f(z)|\leq e^{c_{1}e^{c_{2}d_{o}(z)}}\ \ \ {\rm for\ all}\ \ \ z\in bD\ \! ;
$$
\item[(2)]
for any $z_{1},z_{2}\in bD$ with  $d(z_{1},z_{2})\leq\delta$
$$
|f(z_{1})-f(z_{2})|\leq e^{c_{1}e^{c_{2}\max\{d_{o}(z_{1}),d_{o}(z_{2})\}}}
d(z_{1},z_{2})\ \! .
$$
\end{itemize}
\begin{Th}\label{te1}
There is a constant $c>0$ such that for any
$CR$-function $f$ on $bD$ satisfying conditions (1) and (2) with $c_{2}<c$
there exists $\hat f\in {\cal O}(D)\cap C(\overline{D})$ 
such that
$$
\hat f|_{bD}=f\ \ \ {\rm and}\ \ \ |\hat f(z)|\leq 
e^{\widetilde c_{1}e^{c_{2}d_{o}(z)}}\ \ \ {\rm for\ all}\ \ \ z\in D 
$$
with $\widetilde c_{1}$ depending on $c_{1}$, $c_{2}$, $\delta$ and $c$.
\end{Th}
\begin{R}\label{re1}
{\em (A) If, in addition, $bD$ is smooth of class $C^{k}$, 
$1\leq k\leq\infty$, and $f\in C^{s}(bD)$, $1\leq s\leq k$, then the 
extension $\widehat f$
belongs to ${\cal O}(D)\cap C^{s}(\overline{D})$. This follows from
[HL, Theorem 5.1].\\
(B) Condition (2) means that $f$ is locally Lipschitz with local
Lipschitz constants growing double exponentially.}
\end{R}
\begin{C}\label{c1}
Assume that instead of condition (1) the function $f$ in Theorem \ref{te1} 
satisfies the weaker condition
\begin{itemize}
\item[$(1')$]
$$
|f(z)|\leq e^{\phi(z)}\ \ \ {\rm for\ all}\ \ \ z\in bD
$$
\end{itemize}
where $\phi:M'\to\Re$ is a uniformly continuous function with respect to $d$.

Then there is a constant $C$ (depending on $M,M'$ and $\phi$ only) and
a function $\widehat f\in {\cal O}(D)\cap C(\overline{D})$ such that
$$
\widehat f|_{bD}=f\ \ \ {\rm and}\ \ \
|\widehat f(z)|\leq Ce^{\phi(z)}\ \ \ {\rm for\ all}\ \ \ z\in D
$$
with $C=1$ for $\phi\equiv const$.
\end{C}
\sect{\hspace*{-1em}. Proofs.}
{\bf 2.1. Proof of Theorem \ref{te1}.} Since $r(D)\subset\subset M$, there
is a strictly pseudoconvex domain $S\subset\subset M$ such that
$r(D)\subset\subset S$. Let $S'$ be a connected component of 
$r^{-1}(S)\subset M'$ containing $D$. Then $r:S'\to S$ is an unbranched
covering of $S$. Further, it follows from [Br1, Theorem 2.1] that there
is a function $g\in {\cal O}(S')\cap C(\overline{S'})$ and a constant
$C=C(S',M')$ such that
\begin{equation}\label{eq1}
|g(z)-d_{o}(z)|< C\ \ \ {\rm and}\ \ \ |dg(z)|< C\ \ \ {\rm for\ all}
\ \ \ z\in S'.
\end{equation}
(Here the norm $|\omega(z)|$ of a differential form $\omega$ at $z\in S'$
is determined with respect to the Riemannian metric pulled back from $M$.)
From the first inequality in (\ref{eq1}) one obtains, see [Br1, Example 4.3],
that {\em there is a constant $c>0$ such
that for any $c_{1}>0$ and $0<c_{2}<c$ }
\begin{equation}\label{eq2}
e^{c_{1}e^{c_{2}d_{o}(z)}}\leq\left|e^{c_{1}'e^{c_{2}g(z)}}\right|\leq
e^{c_{1}''e^{c_{2}d_{o}(z)}}\ \ \ {\rm for\ all}\ \ \ z\in S'\ ;
\end{equation}
here $c_{1}', c_{1}''$ are positive constants depending on $c_{1}$, $c_{2}$,
$c$ so that $c_{1}''\to 0$ as $c_{1}\to 0$. We set 
\begin{equation}\label{eq2'}
G_{c_{1},c_{2}}(z):=e^{-c_{1}'e^{c_{2}g(z)}}\ ,\ \ \ z\in S'\ .
\end{equation}

Let us choose $c$ in Theorem \ref{te1} to be the same as in (\ref{eq2}).
Retaining the notation of Theorem \ref{te1} consider the function
$$
f_{1}(z):=f(z)\cdot G_{c_{1}e^{2c_{2}\delta},c_{2}}(z)\ ,\ \ \ z\in bD\ .
$$
\begin{Lm}\label{le1}
$f_{1}$ is a bounded Lipschitz $CR$-function on $bD$.
\end{Lm}
{\bf Proof.} Condition (1) for $f$ and the definition of $f_{1}$ imply that
$$
|f_{1}(z)|\leq 1\ \ \ {\rm for\ all}\ \ \ z\in bD\ .
$$
Thus to show that $f_{1}$ is Lipschitz it suffices to check the 
Lipschitz condition for all pairs $z_{1},z_{2}\in bD$ with 
$d(z_{1},z_{2})\leq\delta$. For such pairs we have
$$
\begin{array}{c}
\displaystyle
|f_{1}(z_{1})-f_{1}(z_{2})|\leq |f(z_{1})-f(z_{2})|\cdot 
|G_{c_{1}e^{2c_{2}\delta},c_{2}}(z_{1})|\ +\\
\\
\displaystyle
|f(z_{2})|\cdot
|G_{c_{1}e^{2c_{2}\delta},c_{2}}(z_{1})-
G_{c_{1}e^{2c_{2}\delta},c_{2}}(z_{2})|
:= I+II\ .
\end{array}
$$
Using condition (2) for $f$, (\ref{eq2}) and the triangle inequality we obtain
\begin{equation}\label{eq3}
I\leq e^{c_{1}(e^{c_{2}\max\{d_{o}(z_{1}),d_{o}(z_{2})\}}-
e^{c_{2}(d_{o}(z_{1})+2\delta)})}d(z_{1},z_{2})\leq d(z_{1},z_{2})\ .
\end{equation}
To estimate $II$ note that according to (\ref{eq1})
and (\ref{eq2}) we have
$$
|dG_{c_{1},c_{2}}(z)|=|-c_{1}'c_{2}G_{c_{1},c_{2}}(z)e^{c_{2}g(z)}dg(z)|
\leq c_{1}'c_{2}Ce^{-c_{1}e^{c_{2}d_{o}(z)}+c_{2}(C+d_{o}(z))}\ .
$$
From here we obtain (for some $c_{3}>0$)
\begin{equation}\label{eq4}
|dG_{c_{1}e^{2c_{2}\delta},c_{2}}(z)|\leq c_{3}e^{-c_{1}
e^{c_{2}(\delta+d_{o}(z))}}
\ \ \ {\rm for\ all}\ \ \  z\in S'\ .
\end{equation}
Further, since $\overline{r(D)}\subset S$ is compact, for a sufficiently 
small $\delta$ the metric $d$ is geodesic in any metric ball 
$B_{\delta}$ on $S'$ of radius $\delta$
centered at a point $D$. (This follows from the definition
of $d$.) Without loss of generality we will assume that this $\delta$ is
the same as in Theorem \ref{te1}. Thus integrating inequality 
(\ref{eq4}) along geodesics in $S'$ we get
\begin{equation}\label{eq5}
II\leq e^{c_{1}e^{c_{2}d_{o}(z_{2})}}c_{3}e^{-c_{1}e^{c_{2}(\delta+
(d_{o}(z_{2})-\delta))}}d(z_{1},z_{2})=c_{3}d(z_{1},z_{2})\ .
\end{equation}
Now the Lipschitz condition for $f_{1}$ follows from inequalities
(\ref{eq3}) and (\ref{eq5}). The fact that $f_{1}$ is $CR$ follows directly
from its definition.\ \ \ \ \ $\Box$

Based on this lemma we reduce the required statement to an extension 
theorem
for the function $f_{1}$. Namely we will show\\
{\bf Claim.} {\em Under the hypotheses of
the lemma there is a function 
$\widehat f_{1}\in {\cal O}(D)\cap C(\overline{D})$ such that }
$$
\widehat f_{1}|_{bD}=f_{1}\ \ \ {\rm and}\ \ \ \sup_{D}|\widehat f_{1}|=
\sup_{bD}|f_{1}|\ .
$$
Then the function 
$\widehat f:=\widehat f_{1}/G_{c_{1}e^{2c_{2}\delta},c_{2}}$ 
is the required extension of Theorem \ref{te1}.

To establish this claim, first, using the McShane theorem [M] let us extend 
$f_{1}$ to a Lipschitz function
$\widetilde f_{1}$ on $S'$ with the same Lipschitz constant as for 
$f_{1}$. Since locally the metric $d$ is equivalent to the Euclidean metric,
by the Rademacher theorem, see, e.g., [Fe, Section 3.1.6],
the function $\widetilde f_{1}$ is differentiable almost everywhere.
In particular, $\overline\partial\widetilde f_{1}$ is a 
$(0,1)$-form on $S'$ whose coefficients in its 
local coordinate representations 
are $L^{\infty}$-functions. Let $\chi_{D}$ be the characteristic function of
$D$. Consider the $(0,1)$-form on $S'$ defined by 
$$
\omega:=\chi_{D}\cdot\overline\partial\widetilde f_{1}\ .
$$
\begin{Lm}\label{le2}
$\omega$ is $\overline\partial$-closed in the distributional sense.
\end{Lm}
{\bf Proof.} We must prove that
\begin{equation}\label{eq6}
\int_{S'}\omega\wedge\overline\partial\phi=0
\end{equation}
for every $(2n-1)$-form $\phi$ of class $C^{\infty}$ with a compact support 
on $S'$ (recall that 
$dim_{\Co}S'=n$). Comparing types of the forms in (\ref{eq6})
we see that, in fact, it suffices to prove the latter identity for
$\phi$ of type $(n,n-2)$. Since this problem is local, it suffices
to prove (\ref{eq6}) for $\phi$ supported in a sufficiently small
neighbourhood of a point of $S'$. Further, by the definition of 
$\omega$ applying the Stokes formula we get that 
identity (\ref{eq6}) is valid for 
$\phi$ supported in a sufficiently
small neighbourhood of a point of $D$. Thus it remains to consider
the case of $\phi$ supported in a sufficiently small neighbourhood 
$U_{z}$ of a boundary point $z\in bD$. (Without loss of generality we may
assume that $U_{z}$ is a coordinate neighbourhood.) Thus we have
$$
\int_{S'}\omega\wedge\overline\partial\phi=
\int_{U_{z}}\omega\wedge\overline\partial\phi=
\int_{U_{z}\cap D}\overline\partial\widetilde f_{1}
\wedge\overline\partial\phi=\int_{U_{z}\cap D}
d(\widetilde f_{1}\cdot\overline\partial\phi)=\int_{U_{z}\cap bD}
f_{1}\cdot\overline\partial\phi=0\ .
$$
We used here that $f_{1}$ is $CR$ and the Stokes formula.
\ \ \ \ \ $\Box$
\begin{R}\label{re2}
{\rm Normally, the Stokes formula is applied to smooth forms.
However, it is also valid for forms with Lipschitz coefficients.
To see this we first apply it to sequences of regularized forms obtained
from these Lipschitz forms and then pass to the limit as the parameter
of the regularization tends to 0. To justify this procedure one uses the
fact that for a Lipschitz function $f$ on
a bounded domain $D\subset\Re^{k}$ the sequence $\{f_{\epsilon}\}$
of regularizations of $f$ converges uniformly to $f$ on every compact
subset of $D$ as $\epsilon\to 0$. 
Moreover, the sequence $\{df_{\epsilon}\}$ is uniformly bounded on every
compact subset of $D$ and converges almost everywhere to $df$ (see, e.g.,
[Fe, Section 4.1.2]).}
\end{R}
\begin{Lm}\label{le3}
There is a bounded continuous function $\widetilde F$ on $S'$
equals 0 on $S'\setminus D$ such that 
$\overline\partial\widetilde F=\omega$ in the distributional sense.
\end{Lm}
{\bf Proof.}
Let us consider a finite open cover ${\cal U}=(U_{i})_{i\in I}$ of a 
neighbourhood $N$ of $\overline{S}$ such that each 
$U_{i}$ is relatively 
compact in a simply connected coordinate chart on $M$ and in these local 
coordinates
is identified with an open Euclidean ball in $\Co^{n}$. By 
${\cal\widetilde U}$ we denote the open cover $(r^{-1}(U_{i}))_{i\in I}$
of $N':=r^{-1}(N)$. In every connected component $V$ of 
$r^{-1}(U_{i})$
we introduce the local coordinates obtained by the pullback of the 
coordinates on $U_{i}$. (Note that $r|_{V}:V\to U_{i}$ is biholomorphic.) 
By the definition of ${\cal\widetilde U}$ in every such $V$ 
the metric $d$ is equivalent to the Euclidean distance on $\Co^{n}$ with 
the constants of the equivalence depending on $U_{i}$ only. 
Since $f_{1}$ is Lipschitz, this implies that in $V$
the form $\omega$ is written as
\begin{equation}\label{eq7}
\omega(z):=\sum_{j=1}^{n}a_{j}(z)d\overline{z}_{j}\ \ \ \ {\rm with}
\ \ \ \ \sup_{z\in V, 1\leq j\leq n}|a_{j}(z)|\leq C\ ;
\end{equation}
here $z_{1},\dots,z_{n}$ are the above introduced local coordinates on $V$
and the constant $C$ is independent of the  choice of $V$ and $U_{i}$. 
Based on
Lemma \ref{le2} and using (\ref{eq7}) one can solve the equation
$\overline\partial F=\omega$ on $V$ to obtain a solution $F_{V}$ which is an 
$L^{\infty}$-function on $V$ satisfying
\begin{equation}\label{eq8}
\sup_{z\in V}|F_{V}(z)|\leq C'
\end{equation}
where $C'$ depends on $C$ and $n$ only, see, e.g., [H, Theorem 6.9]. (Here
$F_{V}$ is the solution in the distributional sense.) Moreover, if we 
identify $V$ with the unit Euclidean ball $B\subset\Co^{n}$ we
can find such a solution $F_{V}$ by the formula
$$
F_{V}(z)=\frac{n!}{(2\pi i)^{n}}\int_{(\xi,\lambda_{0})\in B\times [0,1]}
\omega(\xi)\times
\eta\left((1-\lambda_{0})\frac{\overline\xi-\overline z}{|\xi-z|^{2}}+
\lambda_{0}\frac{\overline\xi}{1-<\overline\xi,z>}\right)\wedge\eta(\xi)
\ ,
$$
see, e.g., [H, Section 4.2];
here for $v=(v_{1},\dots, v_{n})$ and $w=(w_{1},\dots,w_{n})$
$$
\eta(v)=dv_{1}\wedge\cdots\wedge dv_{n}\ ,\ \ \  
<v,w>=\sum_{j=1}^{n}v_{j}\cdot w_{j}\ \ \ {\rm and}\ \ \ |v|^{2}=<v,v>\ .
$$
Since the coefficients of $\omega$ are $L^{\infty}$-functions, 
the above formula implies also that $F_{V}$ is continuous on $V$.
Indeed to show that $F_{V}$ is continuous at $z_{0}\in V$ consider a sequence
$\{z_{i}\}$ convergent to $z$. Without loss of generality we assume
that $\{z_{i}\}$ belongs to the open Euclidean ball $B_{\epsilon}(z_{0})$
centered at $z_{0}$ of radius $\epsilon$ for a sufficiently small
$\epsilon$. Next, we write $F_{V}(z)=F_{1}(z)+F_{2}(z)$ where 
$F_{1}(z)$ is obtained by the integration of the integrand form in the
definition of $F_{V}(z)$ over
$B_{\epsilon}(z_{0})\times [0,1]$ and $F_{2}(z)$ by the integration of this
form over $(B\setminus B_{\epsilon}(z_{0}))\times [0,1]$.
Since the integrand forms for
$F_{2}(z_{i})$
are uniformly bounded on $(B\setminus B_{\epsilon}(z_{0}))\times [0,1]$
and pointwise converge there as $i\to\infty$ to the integrand form for 
$F_{2}(z_{0})$, 
$\lim_{i\to\infty}F_{2}(z_{i})=F_{2}(z_{0})$. To estimate $F_{1}(z_{i})$ we 
use the substitution $w=\xi-z_{i}$ and pass to polar coordinates in the 
obtaining integral. Then it is readily seen that for some $c>0$
$$
|F_{1}(z_{i})|\leq cC\cdot diam(B_{\epsilon}(z_{0}))=2\epsilon cC\ ,
\ \ \
0\leq i<\infty\ .
$$
Therefore $\lim_{i\to\infty}|F_{V}(z)-F_{V}(z_{i})|\leq 4\epsilon cC$ for 
any $\epsilon$, that is $F_{V}$ is continuous at $z_{0}$.

Further, for connected components $V$ and $W$ of $r^{-1}(U_{i})$ and
$r^{-1}(U_{j})$ such that $U_{i}\cap U_{j}\neq\emptyset$
we set 
$$
F_{VW}(z)=F_{V}(z)-F_{W}(z)\ ,\ \ \  z\in V\cap W\ .
$$
Since $\overline\partial F_{VW}=0$ in the distributional sense, 
$F_{VW}\in {\cal O}(V\cap W)$.
Thus considering all possible $V$ and $W$ we get a holomorphic 1-cocycle
$\{F_{VW}\}$ on the cover ${\cal\widetilde U}$ of $N'$. 
Moreover, by (\ref{eq8}) we have
$$
\sup_{V,W, z\in V\cap W}|F_{VW}(z)|\leq 2C'\ .
$$
This implies that the direct image of $\{F_{VW}\}$ with respect to $r$ 
is a holomorphic
1-cocycle on the cover ${\cal U}$ with values in a holomorphic 
Banach vector bundle
with the fibre $l_{\infty}(X)$ where $X$ is the fibre of the 
covering $r: S'\to S$, see the proof of Theorem 2.1 in [Br] for details.
Repeating literally the main argument from the proof of this theorem 
(based on a Banach valued 
version of Cartan's B theorem) together with the fact that
there is a Stein neigbourhood $N_{1}$ of $\overline{S}$ such that
$N_{1}\subset N$ we get holomorphic 
functions $f_{V}\in {\cal O}(V\cap S')$ such that
\begin{itemize}
\item[(1)]
$$
f_{V}(z)-f_{W}(z)=F_{VW}(z)\ ,\ \ \ z\in (V\cap W)\cap S'\ ,\ \ \ {\rm and}
$$
\item[(2)]
$$
\sup_{z\in V\cap S'}|f_{V}(z)|\leq \widetilde C\ 
$$
where $\widetilde C$ depends on $C'$, $N_{1}$ and the cover ${\cal U}$ only.
\end{itemize}

Let us define a continuous function $F$ on $S'$ by the formula
\begin{equation}\label{eq9}
F(z):=F_{V}(z)-f_{V}(z)\ ,\ \ \ z\in V\cap S'\ .
\end{equation}
According to (\ref{eq8}) and condition (2) $F$ is bounded. Also, it satisfies
(in the sense of distributions) the equation $\overline\partial F=\omega$
on $S'$. Since $\omega\equiv 0$ outside $\overline{D}$, the function
$F$ is holomorphic there. Observe that since the boundary of 
$D$ is connected, the set $S'\setminus\overline{D}$ is connected. 
Thus the application of 
Corollary 2.9 of [Br] gives a bounded function $H\in {\cal O}(S')$ such
that $H|_{S'\setminus\overline{D}}=F$. We set $\widetilde F:=F-H$.
Then by the definition $\widetilde F$ is bounded and 
continuous on $S'$ equals 0 on
$S'\setminus D$. Moreover, $\overline\partial\widetilde F=\omega$.
\ \ \ \ \ $\Box$

Using this lemma we define 
$$
\widehat f_{1}(z)=\widetilde f_{1}(z)-\widetilde F(z)\ ,\ \ \ 
z\in\overline{D}\ .
$$
Then 
$$
\widehat f_{1}|_{bD}=f_{1}\ \ \ {\rm and}\ \ \ \
\overline\partial\widehat f_{1}(z)=
\overline\partial\widetilde f_{1}(z)-\overline\partial\widetilde F(z)=
\omega-\omega=0\ \ \ {\rm for}\ \ \ z\in D\ .
$$
This shows that 
$\widehat f_{1}\in {\cal O}(D)\cap C(\overline{D})$. Thus $\widehat f_{1}$ 
is the required holomorphic extension of the function $f_{1}$. To complete
the proof of the Claim it suffices to show that
$$
\sup_{D}|\widehat f_{1}|=\sup_{bD}|f_{1}|\ .
$$
To do this let us consider the product $\widehat f_{1}\cdot G_{c_{1},c_{2}}$ 
where $G_{c_{1},c_{2}}$ is the function from (\ref{eq2'}). Since 
the function $\widetilde f_{1}$ is Lipschitz on $S'$, it satisfies
$$
|\widetilde f_{1}(z)|\leq c_{1}+c_{2}d_{o}(z)\ ,\ \ \ z\in S'\ .
$$
But $\widetilde F$ is bounded on $S'$ and therefore the last inequality
implies that
$$
|\widehat f_{1}(z)|\leq\widetilde c_{1}+c_{2}d_{o}(z)\ ,\ \ \ 
z\in \overline{D}\ .
$$
This and (\ref{eq2}) show that for any $\epsilon>0$ there is a positive
number $R$ such that for any $z\in\overline{D}$ satisfying
$d_{o}(z)\geq R$ one has
$$
|\widehat f_{1}(z)\cdot G_{c_{1},c_{2}}(z)|<\epsilon\ .
$$
In particular, there is an $R_{0}$ such that
$$
\sup_{D}|\widehat f_{1}\cdot G_{c_{1},c_{2}}|=
\sup_{B_{R_{0}}\cap\overline{D}}
|\widehat f_{1}\cdot G_{c_{1},c_{2}}|
$$
where $B_{R_{0}}$ is the open ball on $S'$ centered at
$o$ of radius $R_{0}$. Since $\overline{B}_{R_{0}}\cap\overline{D}$ is
compact and 
$\widehat f_{1}\cdot G_{c_{1},c_{2}}\in {\cal O}(D)\cap C(\overline{D})$, 
the last identity implies that
$$
\sup_{D}|\widehat f_{1}\cdot G_{c_{1},c_{2}}|=
\sup_{\overline{B}_{R_{0}}\cap bD}
|\widehat f_{1}\cdot G_{c_{1},c_{2}}|\leq\sup_{bD}|f_{1}|\ .
$$
Finally, observe that $\widehat f_{1}\cdot G_{c_{1},c_{2}}$ converges 
pointwise to $\widehat f_{1}$ as $c_{1}\to 0$, see (\ref{eq2}). Therefore
we have
$$
|\widehat f_{1}(z)|=
\lim_{c_{1}\to 0}|\widehat f_{1}(z)\cdot G_{c_{1},c_{2}}(z)|
\leq\sup_{bD}|f_{1}|\ ,\ \ \ z\in D\ .
$$
This implies the required identity and completes the proof of the Claim
and therefore of the theorem.
\ \ \ \ \ $\Box$\\
{\bf 2.2. Proof of Corollary \ref{c1}.} 
Let $f$ be a $CR$-function satisfying the hypotheses of Corollary \ref{c1}.
Since $\phi$ is uniformly continuous on $M'$ with respect to the path
metric $d$, there is a constant $C'$ such that
\begin{equation}\label{eq9'}
|\phi(z)|\leq C'd_{o}(z)\ ,\ \ \ z\in M'\ .
\end{equation}
In particular, condition $(1')$ implies condition (1) for $f$. Thus by
Theorem \ref{te1} there exists an extension 
$\widehat f\in {\cal O}(D)\cap C(\overline{D})$ of $f$.

Since $S\subset\subset M$ in the proof of Theorem \ref{te1}
is strictly pseudoconvex, 
it follows from [Br, Theorem 2.1] that {\em
for every function $\phi:S'\to\Re$ uniformly continuous with respect to
the metric $d$ on $S'$ there is a holomorphic function 
$f_{\phi}\in {\cal O}(S')$ and a constant $C=C(\phi,S')$ such that}
$$
|f_{\phi}(z)-\phi(z)|<C\ \ \ {\rm for\ all}\ \ \ z\in S'\ .
$$
Let us consider the function 
$$
\widetilde f(z):=\widehat f(z)\cdot e^{-f_{\phi}(z)}\ ,\ \ \ z\in D\ .
$$
Then by the hypotheses
$$
|\widetilde f(z)|\leq e^{C}\ \ \ {\rm for}\ \ \ z\in bD\ .
$$
From (\ref{eq9'}) and Theorem \ref{te1} we obtain for some $c'>0$ 
(with $c_{2}<c$)
\begin{equation}\label{eq10}
|\widetilde f(z)|\leq e^{c'e^{c_{2}d_{o}(z)}}\ ,\ \ \ z\in D\ .
\end{equation}
Let us take $\widetilde c_{2}$ such that $c_{2}<\widetilde c_{2}<c$ and 
consider the
function $\widetilde f\cdot G_{c_{1},\widetilde c_{2}}$ where 
$G_{c_{1},\widetilde c_{2}}$ is
the function from (\ref{eq2'}). Then from (\ref{eq2}) and (\ref{eq10})
it follows that for any $\epsilon >0$ there is a positive number
$R$ such that for any $z\in\overline{D}$ satisfying $d_{o}(z)\geq R$
one has
$$
|\widetilde f(z)\cdot G_{c_{1},\widetilde c_{2}}(z)|<\epsilon\ .
$$
Now we apply the same argument as at the end of the proof of Theorem
\ref{te1} to get
$$
\sup_{D}|\widetilde f\cdot G_{c_{1},\widetilde c_{2}}|\leq
\sup_{bD}|\widetilde f|\ .
$$
Since $G_{c_{1},\widetilde c_{2}}$ converges pointwise to $1$ as 
$c_{1}\to 0$, from the last inequality we obtain
$$
\sup_{D}|\widetilde f|\leq\sup_{bD}|\widetilde f|\leq e^{C}\ .
$$
From here by the definitions of $\widetilde f$ and $f_{\phi}$ we have
$$
|\widehat f(z)|\leq e^{2C}\cdot e^{\phi(z)}\ ,\ \ \ z\in D\ .
$$
Clearly the above arguments give $C=0$ for $\phi\equiv const$.

This completes the proof of the corollary.\ \ \ \ \ $\Box$
\sect{\hspace*{-1em}. Concluding Remarks.}
Let $r:M'\to M$ be an unbranched covering of a Stein manifold
$M$ with $dim_{\Co}M\geq 2$. We equip $M'$ with a path metric $d$ obtained 
by the pullback of a 
Riemannian metric on $M$. Assume that 
$D\subset\subset M$ is a domain with a connected $C^{k}$-boundary $bD$, 
$1\leq k\leq\infty$. 
We set $D'=r^{-1}(D)$ and $bD'=r^{-1}(bD)$.
Let $\psi :M'\to\Re_{+}$ be a function such that
$\log\psi$ is uniformly continuous with respect to $d$.
For every $x\in M$ we introduce the 
Banach space $l_{p,\psi, x}(M')$, $1\leq p\leq\infty$,
of functions $g$ on $r^{-1}(x)\subset M'$ with norm
$$
|g|_{p,\psi,x}:=\left(\sum_{y\in
r^{-1}(x)}|g(y)|^{p}\psi(y)\right)^{1/p} .
$$
Next, we 
introduce the Banach space ${\cal H}_{p,\psi}(D')$, $1\leq p\leq\infty$, of 
functions $f$ holomorphic on $D'$ with norm
$$
|f|_{p,\psi}^{D}:=\sup_{x\in D}|f|_{p,\psi,x}\ .
$$
In [Br, Theorem 2.7] 
a sharper version of Corollary \ref{c1} for domains $D'$ is proved 
in connection with 
a certain problem posed by Gromov, Henkin and Shubin [GHS]. Namely it was
established that

{\em For every $CR$-function $f\in C^{s}(bD')$, $0\leq s\leq k$, satisfying
$$
f|_{r^{-1}(x)}\in l_{p,\psi,x}(M')\ \ \ {\rm for\ all}\ \ \ x\in D \ \ \
{\rm and}\ \ \
\sup_{x\in bD}|f|_{p,\psi,x}<\infty
$$
there exists a 
function $f'\in {\cal H}_{p,\psi}(D')\cap C^{s}(\overline{D'})$ such 
that $f'|_{bD'}=f$. Moreover, for some $c=c(M',M,\psi,p)$ one has}
$$
|f'|_{p,\psi}^{D}\leq c\sup_{x\in bD}|f|_{p,\psi,x}\ .
$$
It would be interesting to formulate and prove an analog of this result
for other unbounded domains in $M'$. 
A possible formulation of such a result is as follows.

Let $D\subset M'$ be an unbounded domain with a connected smooth boundary
$bD$ such that $r(D)\subset\subset M$. By $dV_{M'}$ and $dV_{bD}$ we
denote the Riemannian 
volume forms on $M'$ and $bD$, respectively, obtained by the pullback of 
a Riemannian metric on $M$. Next, by 
$H_{\psi}^{p}(D)$, $1\leq p\leq\infty$,
we denote the Banach space of holomorphic functions $g$ on $D$ with norm
$$
\left(\int_{z\in D}|g(z)|^{p}dV_{M'}(z)\right)^{1/p}\ .
$$
Also, $L_{\psi}^{p}(bD)$ stands for the Banach space of functions $g$ on 
$bD$ with norm
$$
\left(\int_{z\in D}|g(z)|^{p}dV_{bD}(z)\right)^{1/p}\ .
$$
{\bf Problem.} {\em Let $f\in L_{\psi}^{p}(bD)\cap C(bD)$ be a $CR$-function.
Under what additional conditions on $f$ and $bD$ does there exist
$f'\in H_{\psi}^{p}(D)\cap C(\overline{D})$ such that $f'|_{bD}=f$?}\\

\end{document}